\documentclass[a4paper,11pt]{article}
\usepackage[english]{babel}							
\usepackage{amssymb,amsmath,cancel,mathrsfs} 

\usepackage{epsfig}
\usepackage{slashed}

\def\beq{\arraycolsep1pt\begin{eqnarray*}}
\def\eeq{\end{eqnarray*}}

\newcommand\R[0]{\mathbb{R}}
\newcommand\N[0]{\mathbb{N}}
\newcommand\Z[0]{\mathbb{Z}}
\newcommand\RR[0]{R_0}
\newcommand\ds[0]{\displaystyle}
\newcommand\ee{\varepsilon}
\newcommand\pscal[2]{\left\langle \, #1 \, , \, #2 \, \right\rangle}
\newcommand\cvd[0]{\hfill$\blacksquare$}

\begingroup
\newtheorem{theorem}{Theorem}[section]
\newtheorem{definition}{Definition}[section]
\newtheorem{lemma}[theorem]{Lemma}
\newtheorem{corollary}[theorem]{Corollary}

\endgroup

\title{Periodic impact motions at resonance of a particle bouncing on spheres and cylinders} 

\author{Andrea Sfecci}

\begin{document}

\maketitle

\begin{abstract}
We investigate the existence of periodic trajectories of a particle, subject to a central force, which can hit a sphere, or a cylinder. We will provide also a Landesman-Lazer type of condition in the case of a nonlinearity satisfying a {\em double} resonance condition. Afterwards, we will show how such a result can be adapted to obtain a new result for the impact oscillator at {\em double} resonance.

\end{abstract}

\section{Introduction}

In this paper we are interested in periodic solutions of the differential equation
\begin{equation}\label{eq1}
{\bf x}'' + {\text{f}}(t,|{\bf x}|) {\bf x}=0\,,
\end{equation}
where ${\bf x}\in \R^d$, with $|{\bf x}|\geq \RR$,  and ${\text{f}}:\R\times[\RR,+\infty)\to\R$ is a continuous function, $T$-periodic in the first variable, with $\RR$ a fixed positive constant. We are going to study the existence of {\em bouncing} periodic solutions. In particular we are looking for solutions ${\bf x}:\R\to\R^d$ solving~\eqref{eq1} when $|{\bf x}|>\RR$ and satisfying a {\em perfect bounce condition}. Such a condition describes a {perfectly elastic bounce} on the sphere $\mathbb S^{d-1} =\{{\bf x}\in \R^d \,:\, |{\bf x}|=R_0\}$: the speed has the same value before and after the bounce but the sign of the radial component changes.

\begin{figure}[h]
\centerline{\epsfig{file=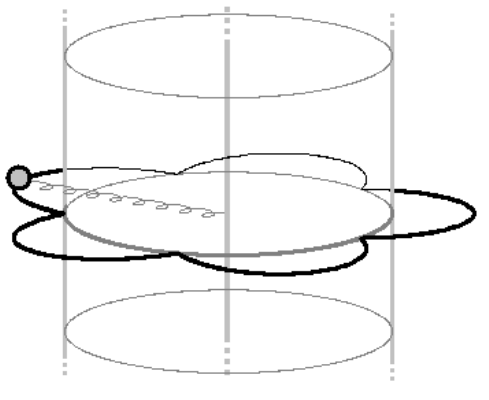, width = 8 cm}}
\caption{a bouncing particle periodically rotating around a cylinder}
\label{cylfig}
\end{figure}

The main results can be applied also to a class of systems with a particle hitting a cylinder $\mathbb S^{d_1-1}\times \R^{d_2}$.
In the case of a proper cylinder (i.e. for $d_1=2$ and $d_2=1$), it models, for example, a particle subject to a periodic central electric field and to an elastic force (see Figure~\ref{cylfig}). Similar situations can be viewed in the case of a proper sphere, for $d=3$.
We will focus our attention, at first, to the case of spheres, postponing the treatment of the case of cylinders to Section~\ref{seccil}.

By the radial symmetry of the equation, every solution of~\eqref{eq1} is contained in a plane, so we can pass to polar coordinates and consider solutions to the following system
\begin{equation}\label{sist1}
\begin{cases}
\ds \rho'' - \frac{L^2}{\rho^3}+f(t,\rho)=0  \qquad \rho>\RR \\[2mm]
\rho^2 \vartheta' = L \,,
\end{cases}
\end{equation}
where $f(t,\rho)=\text{f}(t,\rho)\rho$ and $L\in \R$ is the angular momentum.
The bounce condition could be easily written in the following way
\begin{equation}\label{bounce1}
\rho'(t_0^+)=-\rho'(t_0^-) \text{ if } \rho(t_0)=\RR\,.
\end{equation}
We emphasize that the bounce does not affect the behavior of $\vartheta$.
In this paper we are going to study the existence of rotating periodic solutions, performing a certain number $\nu$ of revolutions around the sphere in the time $kT$ and $T$-periodic in the $\rho$ variable, i.e. such that
\begin{equation}\label{percond}
\begin{array}{l}
\ds \rho(t+T)=\rho(t)\,,\\
\ds \vartheta(t+kT)=\vartheta(t)+2\pi\nu\,. 
\end{array}
\end{equation}

The existence of periodic solutions of large period $kT$ was previously studied by Fonda and Toader in~\cite{FTonProc} in the setting of a Keplerian-like system where the planet is viewed as a point (see also~\cite{FTonJDE,FTZ} for other situations).

Problems modeling the motion of a particle hitting some surfaces have been widely studied in literature in different situations, see e.g.~\cite{BattF,BP2,CP2,J2,Lamba,LM,RHT}. The simpler system with impacts is given by the so-called impact oscillator (see e.g.~\cite{Bap,BF,J,QT2}) where a particle hits a wall attracted towards it by an elastic force. The existence of bouncing periodic solutions of such systems has been discussed e.g. in~\cite{BF,FS3,LM,O,OonPM,Q,QT2,SQ}. However, to the best of our knowledge, it seems that similar existence results on rotating periodic solutions with impact on spheres (or cylinders) of positive radius have not been presented yet.

\medbreak

Let us now explain in details what we mean by the term ``bouncing solution'', borrowing the definition given by Bonheure and Fabry in~\cite{BF}, which we recall for the reader convenience. In it, the constant $w$ indicates the $x$-coordinate of the {\em wall} against which the solution bounces.

\begin{definition}\label{bsol}
Consider a scalar second order differential equation
$$
x'' + p(t,x)=0 \,,
$$
where $p:\R\times[w,+\infty)\to\R$ is a continuous function. A $w$-{\em bouncing solution} is a continuous function $x(t)$, defined on some interval $(a,b)$, such that $x(t)\geq w$ for every $t\in(a,b)$, satisfying the following properties:
\begin{itemize}
	\item[\textit{i.}] if $t_0\in (a,b)$ is such that $x(t_0)>w$, then $x(t)$ is twice differentiable at $t=t_0$, and $x''(t_0)+ p(t_0,x(t_0))=0$;
	\item[\textit{ii.}] if $t_0\in (a,b)$ is such that $x(t_0)=w$ and, in a neighborhood of $t_0$, $x(t)>w$ for $t\neq t_0$, then $x'(t_0^-)$ and $x'(t_0^+)$ exist and $x'(t_0^-)=-x'(t_0^+)$;
	\item[\textit{iii.}] if $t_0\in (a,b)$ is such that $x(t_0)=w$ and, either $x'(t_0^-)$, or $x'(t_0^+)$, exists and is different from $0$, then, in a neighborhood of $t_0$, $x(t)>w$ for $t\neq t_0$.
	\item[\textit{iv.}] if $x(t)=w$ for all $t$ in a non-trivial interval $I\subseteq(a,b)$, then $p(t,w)\ge0$ for every $t\in I$.
\end{itemize}
\end{definition}

\medskip

Treating a Keplerian system like~\eqref{sist1}, we will say that $(\rho,\vartheta)$ is a $w$-bouncing solution if $\rho$ is a $w$-bouncing solution of the first differential equation.

We are now ready to state one of the main results of this paper.

\begin{theorem}\label{main}
Assume
\begin{equation}\label{asympt}
\check \mu \leq \liminf_{\rho\to+\infty} \frac{f(t,\rho)}{\rho} \leq \limsup_{\rho\to+\infty} \frac{f(t,\rho)}{\rho} \leq \hat\mu \,,
\end{equation}
with
\begin{equation}\label{nonres}
\left( \frac{N \pi}{T} \right)^2 < \check\mu \leq \hat\mu < \left( \frac{(N+1) \pi}{T} \right)^2 \,,
\end{equation}
for a suitable integer $N$. Then, for every integer $\nu>0$, there exists an integer $k_\nu>0$, such that for every integer $k\geq k_\nu$ there exists at least one periodic $R_0$-bouncing solution $(\rho,\vartheta)$ of~\eqref{sist1} with period $kT$, which makes exactly $\nu$ revolutions around the origin in the period time $kT$, i.e. satisfying~\eqref{percond}.
\end{theorem}

The statement of the theorem requires that the nonlinearity $f$ has a {\em nonresonant}  asymptotically linear growth at infinity in the following sense. The constants $\mu_j=(\pi j / T)^2$ in~\eqref{nonres} are the values of the vertical asymptotes $\mu=\mu_j$ of the $j$-th curve of the periodic Dancer-Fu\v c\'{\i}k spectrum associated to the asymmetric oscillator $x'' + \mu x^+ - \nu x^- =0$. The first equation in~\eqref{sist1} presents a {\em wall} at $\rho=R_0$ against which the particle bounces. Such a wall can be approximatively modeled as a spring with a very large elasticity constant (see e.g.~\cite{BF,BP2,LM,SQ}).
For this reason, it is natural to require that the nonlinearity $f$ satisfies an asymptotic behavior at infinity as in~\eqref{asympt}. Conditions like~\eqref{asympt} have already been introduced treating scalar equations with a singularity, see e.g.~\cite{DMM,FG2,FTonNA}.

\medbreak

We will see also, in Theorem~\ref{LLtheorem}, how we can relax condition~\eqref{nonres} introducing a Landesman-Lazer type of condition thus obtaining a similar existence result for nonlinearities {\em next to} resonance. In particular we will introduce a double resonance condition for nonlinearities satisfying~\eqref{asympt} with $\check\mu=\mu_N$ and $\hat\mu=\mu_{N+1}$. Double Landesman-Lazer type of condition has been treated in other situations, see e.g.~\cite{Fabry,FaFo,FG}.

Moreover, following the proof of Theorem~\ref{LLtheorem}, we will obtain Corollary~\ref{LLcorollary} which extends to the resonant case a previous result obtained by Fonda and the author in~\cite{FS3} for impact oscillators.

\medbreak

The next section is devoted to the proof of Theorem~\ref{main}. In Section~\ref{LLsec} we will present  the results {\em next to} resonance and we will show how to modify the proof of Theorem~\ref{main} in order to prove such a result.
Then, in Section~\ref{seccil}, we will present how to extend the applications to cylinders.

\section{Nonresonant case, proof of Theorem~\ref{main}}

Let us consider the change of coordinate $\rho=r+\RR$. System~\eqref{sist1} is equivalent to
\begin{equation}\label{sist}
\begin{cases}
\ds  r'' - \frac{L^2}{(r+\RR)^3}+f(t,r+\RR)=0  \qquad r>0\\
\ds \vartheta' = \frac{L}{(r+\RR)^2} \,.
\end{cases}
\end{equation}
Let us define the function $g:\R\times\R\times[0,+\infty)\to\R$ as
\begin{equation}\label{whoisg}
g(L,t,r)= - \frac{L^2}{(r+\RR)^3}+f(t,r+\RR)\,,
\end{equation}
so that the first differential equation in~\eqref{sist} becomes
\begin{equation}\label{BP}
r'' + g(L,t,r)=0
\end{equation}
and it is easy to verify that
\begin{equation}\label{asymptG}
\check\mu \leq \liminf_{r\to+\infty} \frac{g(L,t,r)}{r} \leq \limsup_{r\to+\infty} \frac{g(L,t,r)}{r} \leq \hat\mu
\end{equation}
uniformly in $t\in [0,T]$ and $L$ in a compact set.
In what follows we will assume that $L$ varies in a compact set containing zero, but it is not restrictive to assume $L$ to be non-negative. So, fixing $L_0>0$, in what follows, we will always assume $L\in[0,L_0]$.

We are looking for $0$-bouncing solutions to~\eqref{BP}
such that

$$
\begin{array}{l}
\ds r(t+T)=r(t)\,,\\
\ds \vartheta(t+kT)=\vartheta(t)+2\pi\nu\,, \quad  \text{ with } k,\nu \in \Z\,.
\end{array}
$$

We define, for a fixed small $\delta>0$, for every $n\in\N$, the functions $g_n: [0,L_0]\times\R\times\R\to\R$
\begin{equation}\label{gn}
g_n(L,t,x)=
\begin{cases}
g(L,t,x) & x\geq 1/n\\
nx \big( g(L,t,x) + \delta \big) - \delta & 0<x<1/n\\
nx-\delta & x\leq 0 \,
\end{cases}
\end{equation}
and we consider the differential equations
\begin{equation}\label{approx}
x'' + g_n(L,t,x)=0\,,
\end{equation}
where, here, $x$ varies in $\R$.
Let us spend few words to motivate the introduction of the small constant $\delta>0$, which can appear unessential: it will be useful to simplify the proof of the validity of the fourth property, which appears in Definition~\ref{bsol}, for the bouncing solution we are going to find (cf.~\cite{FS3}).

It is well known that such a differential equation has at least one $T$-periodic solution when $n$ is large enough and $L$ is fixed (cf.~\cite{DI,FH,FS1}). In the next section, we will prove the existence of a common a priori bound.

Let us introduce the set of $C^1$ functions which are $T$-periodic
$$
C^1_P = \{ x\in C^1([0,T]) \,:\, x(0)=x(T) \,, x'(0)=x'(T) \}\,.
$$
We are going to look for an open bounded set $\Omega\subset C^1_P$, containing $0$, such that any periodic solution of~\eqref{approx} belongs to $\Omega$.

Then, we will prove in Section~\ref{degree} the existence of periodic solutions of~\eqref{BP} by a limit procedure: we will show that, for every integer $\nu>0$ and for every integer $k$ sufficiently large, there exists a sequence $\big(L^{k,\nu}_n\big)_n\subset[0,L_0]$ and a sequence of solutions $\big(x^{k,\nu}_n\big)_n$ of~\eqref{approx}, with $L=L_n$, converging respectively to $L^{k,\nu}$ and $x^{k,\nu}$, where $x^{k,\nu}$ is the desired bouncing solution. The proof of Theorem~\ref{main} will be concluded easily in Section~\ref{conclusion}.

\subsection{The a priori bound}\label{secnonres}

The proof makes use of some phase-plane techniques, so it will be useful to define the following subsets of $\R^2$:
$$
\Pi^-=\{ (x,y)\in\R^2 \,:\, x \leq 0\}
\,,\quad
\Pi^+=\{ (x,y)\in\R^2 \,:\, x \geq 0\}
$$
and the open balls centered at the origin
$$
B_s = \{ (x,y)\in\R^2 \,:\, x^2+y^2 < s^2\}\,.
$$

Let $x\in C^1_P$ be a solution of~\eqref{approx} such that there exist some instants $t_1,t_2,t_3\in[0,T]$ such that, $x(t_i)=0$ and $(-1)^i x'(t_i)>0$. Assume moreover that $x(t)<0$ for every $t\in(t_1,t_2)$ and $x(t)>0$ for every $t\in(t_2,t_3)$. Let us first consider the interval $[t_1,t_2]$. There exists a positive constant $c$ such that the orbit $(x(t),x'(t))$ in the phase-plane, in this interval, consists of a branch of the ellipse
\begin{equation}\label{ellipse}
y^2 + n x^2 - 2\delta x = c^2 \,.
\end{equation}
In particular $x'(t_1)=-c$ and $x'(t_2)=c$, and for every $t\in[t_1,t_2]$ one has $x(t)\geq
\frac{\delta-\sqrt{\delta^2+c^2 n}}{n}>-\frac{c}{\sqrt n}$. A computation shows that, for every $c>0$,
$$
t_2-t_1 = \frac{1}{\sqrt{n}} \left[ \pi - 2 \arcsin \left(\frac{\delta}{\sqrt{nc^2 +\delta^2}} \right) \right] < \frac \pi{\sqrt n}\,.
$$

Call $\Sigma^n_{c}$ the open region delimited by the ellipse in~\eqref{ellipse}. Besides, it is possible to find $n_0=n_0(c)$ large enough to have
\begin{equation}\label{nzero}
\frac{1}{n}\left(\delta-\sqrt{\delta^2+c^2 \,n}\,\right)>-R_0/2 \,, \quad \text{ for every } n\geq n_0(c)\,.
\end{equation}
It is useful for our purposes to define the following set
\begin{equation}\label{goodchoice}
\Xi_c^{n} = (\Pi^-\cap \Sigma^{n}_{c}) \cup (\Pi^+\cap B_{c})\,.
\end{equation}
In particular we have $\Xi_c^{n+1} \subset \Xi_c^n$ and
$$
\Xi_c^{n_0(c)} \subset (-R_0/2,c )\times (-c,c)\,.
$$
We define the following set of periodic functions
\begin{equation}\label{omeganc}
\Omega_c^n=\{x\in C^1_P \,:\, (x(t),x'(t))\in\Xi_c^n \text{ for every } t\in[0,T] \} \,.
\end{equation}

Let us now focus our attention on the second interval $[t_2,t_3]$.
It is possible to verify that there exists $\bar\chi>0$ (which can be chosen independently of $L$ and $n$) sufficiently large with the following property: if $x$ satisfies $(x'(t))^2 + (x(t))^2 > \bar\chi^2$ for every $t\in[t_2,t_3]$, then (cf.~\cite{FM})
$$
\frac{\pi}{\sqrt{\hat\mu}} \leq t_3-t_2 \leq \frac{\pi}{\sqrt{\check\mu}}\,.
$$

Hence, we can conclude that any solution to~\eqref{approx}, such that $x^2+x'\;^2>\bar\chi^2$ when $x\geq 0$, must rotate in the phase plane, and it needs a time $\tau$, such that
\begin{equation}\label{tau}
\tau \in \left[ \frac{\pi}{\sqrt{\hat\mu}} \,, \frac{\pi}{\sqrt{\check\mu}} + \frac{\pi}{\sqrt{n}}\,\right] \,,
\end{equation}
to complete a rotation. So, by~\eqref{nonres}, choosing $n$ large enough we can find that it performs more than $N$ rotations and less than $N+1$ rotations in the phase plane. In particular it cannot perform an integer number of rotations.

Moreover, by assumption~\eqref{asympt}, if we consider a solution $x$ to~\eqref{approx},  satisfying $x^2 + x'^2>\bar\chi^2$ when $x\geq 0$,
introducing polar coordinates $(\varrho,\theta)$, we can find (enlarging $\bar\chi$, if necessary) the following uniform bound for the radial and angular velocity when $x$ is positive:
$$
-\theta'(t) = \frac{x'(t)^2 + x(t)g_n(L,t,x(t))}{x(t)^2+x'(t)^2} > \theta_0 > 0 \,,
$$
$$
|\varrho'(t)|=\frac{|x'(t)(x(t)-g_n(L,t,x(t)))|}{\sqrt{x(t)^2+x'(t)^2}} < l_0\, \varrho(t) \,.
$$
thus finding
\begin{equation}\label{rt}
\left| \frac{d\varrho}{d(-\theta)} \right| < \frac{l_0}{\theta_0} \varrho = K\varrho \,.
\end{equation}

\medbreak

Let us now state a lemma which will be useful also in Section~\ref{LLsec}.

\begin{lemma}\label{lemmaspirale}
For every $\chi\geq \bar\chi$ there exists $R=R(\chi)>0$ and an integer $n_1=n_1(\chi)$ such that every solution $x$ to \eqref{approx}, with $n>n_1$ satisfying $(x(t_0),x'(t_0))\in \overline B_\chi^+ = \overline B_\chi \cap \Pi^+$ at a certain time $t_0\in[0,T]$, is such that $(x(t),x'(t))\in \Xi_R^{n_1}$ for every $t\in[t_0,t_0+T]$.
\end{lemma}

\textit{Proof.}
It is not restrictive to consider a solution such that $(x(t_0),x'(t_0))\!\in \partial B_{\chi} \cap \Pi^+$ and  $(x(s),x'(s)) \notin B_\chi \cap \Pi^+$, for every $s\in(t_0,t_0+T]$. Such a solution rotates clockwise in the phase-plane and will vanish for the first time at $t_1>t_0$, thus having $(x(t_1),x'(t_1))=(0,-y_1)$ such that, by~\eqref{rt}, $0<y_1< e^{K\pi} \chi$. The solution vanishes again at $t_2>t_1$ such that  $(x(t_2),x'(t_2))=(0,y_1)$.
Then, the solution will perform a complete rotation in the interval $[t_2,t_3]$ thus obtaining
$x'(t_3) < e^{K\pi} x'(t_2)$. In the time interval $[t_0,t_0+T]$ the solution cannot perform more than $N+1$ rotations, so choosing $R(\chi) =  e^{(N+2)K\pi} \chi$ and $n_1(\chi) = n_0(R(\chi))$, as in~\eqref{nzero}, we conclude the proof of the lemma.
 \cvd

\medbreak

With a similar reasoning, we can prove that if $x'$ vanishes at a time $\tau_0$ with $x(\tau_0)=\bar x>0$ then the solution $x$ will vanish for the first time at $\tau_1$ with $-x'(\tau_1) < \bar x e^{K\pi/2}$. Then, the solution will reach a negative minimum at $\tau_2$ and again will vanish at $\tau_3$ with $x'(\tau_1) = -x'(\tau_1)$. Recalling that, in $\Pi^-$, the orbit of $x$ is contained in a certain ellipse of equation~\eqref{ellipse} we have $x(\tau_2)>-x'(\tau_3)/\sqrt{n}$. Hence, we have immediately, setting $C=e^{K\pi/2}$
\begin{equation}\label{goodderivative}
\frac1C \|x\|_\infty \leq \| x'\|_\infty \leq C \|x\|_\infty \quad\text{and}\quad x(t)>-\frac{\|x'\|_\infty}{\sqrt n} \geq - \frac{C \|x\|_\infty}{\sqrt n} \,.
\end{equation}

Such estimates will be useful in Section~\ref{LLsec}.

\medbreak

The set $\Xi=\Xi_R^{n_1}$ provided by the previous lemma is the a priori bound we were looking for. In fact, suppose to have a $T$-periodic solution such that $(x(\zeta),x'(\zeta))\notin\Xi$ at a certain time $\zeta\in[0,T]$.
If the solution remains outside $B_{\bar\chi}^+$ then it cannot perform an integer number of rotations around the origin in the period time $T$. Hence, the solution must enter the set at a certain time $\zeta'>\zeta$, and the previous lemma gives us a contradiction.

Summing up, in this section, we have proved the following estimate.

\begin{lemma}\label{finallemma}
There exists a positive integer $\bar n$ and an open bounded set $\Xi\subset (-R_0/2,+\infty)\times\R$ such that every $T$-periodic solutions to~\eqref{approx} with $n>\bar n$ belong to
$$
\Omega = \{ x\in C^1_P \,:\, (x(t),x'(t))\in \Xi \text{ for every } t\in[0,T]\,\}\,.
$$
\end{lemma}

\subsection{Degree theory}\label{degree}

It is well-known that the existence of periodic solutions of equation~\eqref{approx} is strictly related to the existence of a fixed point of a completely continuous operator $\Psi_{L,n}: C^1_P \to C^1_P$ (see e.g.~\cite{FTonProc}),
$$
\Psi_{L,n}= (\mathcal L- \sigma I)^{-1} ( \mathcal N_{L,n} -\sigma I )
$$
where $\mathcal L: D(\mathcal L) \to L^1(0,T)$ is defined in
$D(\mathcal L) = \{ x\in W^{2,1}(0,T) \,:\, x(0)=x(T) \,, x'(0)=x'(T)\}$
 as $\mathcal L x = x''$ and $\sigma$ does not belong to its spectrum, $(\mathcal N_{L,n} x) (t) = - g_n(L,t,x(t))$ is the so called Nemytzkii operator and $I$ is the identity operator.

By classical results (see e.g.~\cite{DI,FH}) one has
\begin{equation}\label{deg}
d_{LS}(I-\Psi_{L,n},\Omega)\neq 0 
\end{equation}
for every $n\geq n_0$, with $\Omega$ given by Lemma~\ref{finallemma}.

Using the continuation principle one has, for every $n\geq n_0$, that there exists a continuum $\mathcal C_n$ in $[0,L_0]\times\Omega$, connecting $\{0\}\times\Omega$ to $\{L_0\}\times\Omega$, whose elements $(L_o,x_o)$ are such that $x_o$ is a solution of $x_o''+g_n(L_o,t,x_o)=0$  (see~\cite{FTonProc} for a similar approach). The function $\Theta: [0,L_0]\times\Omega \to \R$
$$
\Theta(L,x)=\int_0^T \frac{L}{(R_0+x(t))^2} \, dt
$$
is well defined
 and continuous, being $x(t)>-R_0/2$, by Lemma~\ref{finallemma}. In particular, one has $\Theta(0,x)=0$ and 
\begin{equation}\label{thetaawayfromzero}
\frac{T}{(R_0+R)^2}\, L <\Theta(L,x)< \frac{4T}{R_0^2} \, L\,,
\end{equation}
where $R$ is the constant provided by Lemma~\ref{lemmaspirale}.
Hence, for every integer $\nu>0$, there exists $k_\nu$ with the following property: for every $k\geq k_\nu$ and for every $n\geq n_0$, there exists $(L^{k,\nu}_n, x^{k,\nu}_n)\in \mathcal C_n$ with $\Theta(L^{k,\nu}_n,x^{k,\nu}_n)=2\pi\nu/k$.

Fix now $\nu$ and $k\geq k_\nu$ and consider the sequences $(L^{k,\nu}_n)_n$ and $(x^{k,\nu}_n)_n$. Let us simply denote them by $(L_n)_n$ and $(x_n)_n$. Both the sequences are contained in a compact set respectively of $\R$ and $C^1_P$, so there exist, up to subsequences, $\overline L>0$ and $\overline x$ such that $L_n\to \overline L$ and $x_n\to \overline x$ uniformly (the estimate in~\eqref{thetaawayfromzero} gives us that $\overline L$ is positive). Moreover, by continuity of $\Theta$, we have $\Theta(\overline L, \overline x)=2\pi\nu/k$. We have to prove that $\overline x$ is a bouncing solution of the differential equation~\eqref{BP}, where $L=\overline L$. The proof follows the same procedure of the one in~\cite[p.~185--188]{FS3} considering the sequence of approximating differential equations
$$
x'' + f_n(t,x)=0 \qquad\text{with } f_n(t,x)=g_n(L_n,t,x)\,.
$$
For briefness we refer to that paper for the proof of this part.

\subsection{Conclusion}\label{conclusion}

We have found in the previous section, for every integer $\nu$ and for every integer $k$ sufficiently large, a solution $x$ of~\eqref{BP} and $L\in(0,L_0]$ such that $\Theta(L,x)=2\pi\nu/k$. So, defining for a certain $\theta_0\in[0,2\pi)$,
$$
\theta(t)=\theta_0 + \int_0^t \frac{L}{(R_0+x(t))^2} \, dt
$$
we have that $(x,\theta)$ is a bouncing periodic solution of~\eqref{sist} and using the change of coordinate $\rho=r+R_0$ we find the bouncing  solution of~\eqref{sist1} satisfying~\eqref{percond}. The proof of Theorem~\ref{main} is thus completed.

\section{Nonlinearities next to resonance,\\ a double Landesman-Lazer type of condition}\label{LLsec}

In this section we will see how we can relax the hypotheses of Theorem~\ref{main} in order to treat the situation of a nonlinearity which has an asymptotically linear growth next to resonance. 
We are going to provide a Landesman-Lazer type of condition for the case when the nonlinearity $f$ satisfies~\eqref{asympt} with $\check\mu=\mu_{N}$ and $\hat\mu=\mu_{N+1}$, where $\mu=\mu_{j}$ is the vertical asymptote of the $j$-th curve of the Dancer-Fu\v c\'{\i}k spectrum. Such a situation has been often called as a {\em double resonance} (see e.g.~\cite{Fabry,FaFo,FG}). 
A {\em one-side} Landesman-Lazer condition for the scalar differential equation with singularity has been provided by Fonda and Garrione in~\cite{FG2}.
This section has been inspired by such a paper and some steps of the proof of Lemma~\ref{lemmaLL} will appear similar. The main novelty occurs in the estimate in~\eqref{cr1} and its proof, permitting us to treat a double resonance situation. In particular, the validity of~\eqref{cr1}, permits us to obtain a Landesman-Lazer condition involving the function
\begin{equation}\label{psigood}
\psi_j(t) =
\sin(\sqrt{\mu_{j}} \, t) \quad \text{ with } t\in \left[0,\frac{T}{j}\right]
\end{equation}
extended by periodicity to the whole real line, while, in~\cite{FG2}, the Landesman-Lazer condition is weaker, being related to a function of the type
\begin{equation}\label{psibad}
\tilde \psi_j(t) =
\begin{cases}
\sin(\sqrt{\mu_{j}} \, t) & t\in \left[0,\frac{T}{j}\right] \\[2mm]
0 & t\in \left[\frac{T}{j}, T\right]
\end{cases}
\end{equation}
extended by periodicity.

We will prove in this section the following result.

\begin{theorem}\label{LLtheorem}
Assume that 
there exists a constant $\hat\eta$ such that, for $N>0$,
$$
\mu_N x -\hat\eta \leq f(t,x)\leq \mu_{N+1} x +\hat\eta
$$
for every $t\in[0,T]$ and every $x>R_0$. Moreover, for every $\tau\in[0,T]$,
\begin{equation}\label{LLcond}
\int_0^T \limsup_{x\to +\infty} ( f(t,x+R_0) - \mu_{N+1} x) \psi_{N+1}(t+\tau) \, dt <0\,
\end{equation}
and
\begin{equation}\label{LLcond2}
\int_0^T \liminf_{x\to +\infty} ( f(t,x+R_0) - \mu_N x ) \psi_{N}(t+\tau) \, dt > 0\,.
\end{equation}
Then, for every integer $\nu>0$, there exists an integer $k_\nu>0$, such that for every integer $k\geq k_\nu$ there exists at least one periodic $R_0$-bouncing solution $(\rho,\vartheta)$ of~\eqref{sist1} with period $kT$, which makes exactly $\nu$ revolutions around the origin in the period time $kT$, i.e. satisfying~\eqref{percond}.
\end{theorem}

\noindent The case $N=0$ will be briefly treated at the end of this section.

First of all we need to introduce a sequence of approximating equations with non-resonant nonlinearities. Then, we will look for a common a priori bound.

It is possible to find, when $n$ is chosen large enough, a constant $\kappa_n>1$ such that
$$
\mu_N < n \, \frac{(\sqrt{\kappa_n}-1)^2}{\kappa_n} < \mu_{N+1}\,.
$$
In this way, we have
\begin{equation}\label{kappan}
\begin{array}{rcl}
\ds \frac{T}{(N+1)\pi} = \frac{1}{\sqrt{\mu_{N+1}}} &<& \ds
\frac{1}{\sqrt{\kappa_n \mu_{N+1}}} + \frac{1}{\sqrt{n}}\\
&<& \ds \frac{1}{\sqrt{\kappa_n \mu_{N}}} + \frac{1}{\sqrt{n}}
< \frac{1}{\sqrt{\mu_{N}}} = \frac{T}{N \pi}\,.
\end{array}
\end{equation}

We define the function $g_n$, similarly as in~\eqref{gn},
\begin{equation}\label{gn2}
g_n(L,t,x)=
\begin{cases}
\tilde g_n(L,t,x) & x\geq 1/n\\
nx \big( \tilde g_n(L,t,x) + \delta \big) - \delta & 0<x<1/n\\
nx-\delta & x\leq 0 \,
\end{cases}
\end{equation}
where
$$
\tilde g_n(L,t,x)= - \frac{L^2}{(x+\RR)^3}+\kappa_n f(t,x+R_0)\,.
$$
We will consider the differential equation
\begin{equation}\label{approx2}
x'' + g_n(L,t,x)=0.
\end{equation}

Notice that
\begin{equation}\label{asympt2}
\kappa_n \mu_N \leq \liminf_{\xi\to+\infty} \frac{g_n(L,t,\xi)}{\xi} \leq \limsup_{\xi\to+\infty} \frac{g_n(L,t,\xi)}{\xi} \leq \kappa_n \mu_{N+1} \,,
\end{equation}
thus giving us, by the previous computation in~\eqref{kappan}, that the nonlinearities $g_n$ are nonresonant. So, if we find a common a priori bound $\Omega$ as in Lemma~\ref{finallemma}, uniform in $L$ and $n$, for every $T$-periodic solution of~\eqref{approx2}, then~\eqref{deg} holds, thus permitting us to end the proof of the theorem as in Section~\ref{degree} and~\ref{conclusion}.
The Landesman-Lazer type of conditions introduced in~\eqref{LLcond} and~\eqref{LLcond2} are needed in order to find the common a priori bound for every $n$ sufficiently large.

We can assume without loss of generality that the constant $\hat\eta$ in the statement of the theorem guarantees also that, when $x> R_0$, for every $L$, $t$ and $n$,
\begin{equation}\label{etacontrol}
g_n(L,t,x) \leq \kappa_n \mu_{N+1} x + \hat\eta \,.
\end{equation}

The reader will notice that the next result comes out free by the proof of the previous theorem (setting $L=0$ everywhere). Such a corollary extends a previous result provided by Fonda and the author in~\cite{FS3}.

\begin{corollary}\label{LLcorollary}
Assume
that there exists a constant $\hat\eta$ such that
$$
\mu_N x - \hat\eta \leq g(t,x)\leq \mu_{N+1} x +\hat\eta
$$
for every $t\in[0,T]$ and every $x>0$. Moreover, for every $\tau\in[0,T]$,
$$
\int_0^T \limsup_{x\to +\infty} ( g(t,x) - \mu_{N+1} x) \psi_{N+1}(t+\tau) \, dt <0\,,
$$
and
$$
\int_0^T \liminf_{x\to +\infty} ( g(t,x) - \mu_N x) \psi_N(t+\tau) \, dt >0\,.
$$
Then, there exists a $0$-bouncing solution for equation $x''+g(t,x)=0$.
\end{corollary}

\medbreak

Repeating the reasoning explained in Section~\ref{secnonres}, providing the estimate in~\eqref{tau}, we can find, for every $\varepsilon>0$, a value $\chi_\varepsilon>0$ such that every solution $x$ of~\eqref{approx2}, satisfying $x^2 + x' \, ^2 > \chi_\varepsilon^2$ when $x\geq 0$, must rotate in the phase plane spending a time $\tau= (t_3-t_2)+(t_2-t_1)$ with
\begin{equation}\label{goodtime}
t_3 - t_2 \in \left( \frac{\pi}{\sqrt{\mu_{N+1}}} -\varepsilon \,,\, \frac{\pi}{\sqrt{\mu_N}} + \varepsilon \right) .
\end{equation}
When $\varepsilon$ is chosen sufficiently small, one has that a $T$-periodic solution $x$ of~\eqref{approx2}, such that $x^2 + x' \, ^2 > \chi_\varepsilon^2$ when $x\geq 0$, must perform exactly $N$ or $N+1$ rotations around the origin. The previous reasoning holds uniformly for every $n$ and every $L$, so we can fix such a suitable $\varepsilon$ and find the constant $\bar \chi=\chi_\varepsilon$.

We underline that Lemma~\ref{lemmaspirale} still holds under the hypotheses of Theorem~\ref{LLtheorem}, too. We underline that~\eqref{goodderivative} remains valid.

\medbreak

The needed a priori bound is given by the following lemma.

\begin{lemma}\label{lemmaLL}
There exist $\bar R\geq R(\bar\chi)$ (given by Lemma~\ref{lemmaspirale}) and $\bar n\geq n_0(\bar R)$, as in~\eqref{nzero}, such that, for every $n>\bar n$ and every $L\in[0,L_0]$, any $T$-periodic solution $x$ of~\eqref{approx2} is such that $(x(t),x'(t)) \in  B_{\bar R}^+$ when $x(t)\geq0$. In particular, we immediately have $x\in\Omega_{\bar R}^{\bar n}$. 
\end{lemma}

\medbreak

\textit{Proof.} 
Suppose by contradiction that there exist an increasing sequence $(R_m)_m$ with $R_m>R(\bar\chi)$ and $\lim_m R_m = +\infty$, an increasing sequence $(n_m)_m$ of integers $n_m>n_0(R_m)$, a sequence $(L_m)_m\subset[0,L_0]$, a sequence $(x_m)_m$ of solutions to
\begin{equation}\label{gnm}
x_m'' + g_{n_m}(L_m,t,x_m)=0
\end{equation}
and a sequence of times $(t^1_m)_m\subset[0,T]$ such that $(x(t_m^1),x'(t_m^1))\notin B_{R_m}$ and $x(t_m^1)\geq0$.

Being $R_m > R(\bar\chi)$, thanks to Lemma~\ref{lemmaspirale}, all the solutions cannot enter $B_{\bar\chi}^+$, so that they must perform exactly $N$ or $N+1$ rotations around the origin.

We define the sequence of functions
$$
v_m = \frac{x_m}{\|x_m\|_\infty}
$$
which are solutions of
\begin{equation}\label{seq}
v_m(t)'' + \frac{g_{n_m}(L_m,t,x_m(t))}{\|x_m\|_\infty}=0\,.
\end{equation}
By~\eqref{goodderivative}, we have
\begin{equation}\label{goodvm}
-\frac{1}{\sqrt{n_m}} \leq v_m(t)\leq 1\, \text{ for every } t\in[0,T] \text{ and } \frac 1{C_0} < \|v_m'\|_\infty < C_0 \,,
\end{equation}
for a suitable $C_0>0$. We have immediately
that $(v_m)_m$ is bounded in $H^1(0,T)$ so that, up to subsequences we have
$v_m\to v$ weakly in $H^1$ and uniformly, moreover we can assume that $L_m\to L_\dag$. In particular $v\neq 0$, being $\|v\|_\infty =1$, and it is non-negative and $T$-periodic.

\medbreak

We assume that, up to subsequence, all the solutions make exactly $N+1$ rotations. We will discuss the other situation later on.

We are going now to prove that $v$ solves $v'' + \mu_{N+1} v =0$ for almost every $t$. Call
$$
\mathcal I^+ = \{ t\in\R : v(t)>0 \}\,,
$$
which is an at most countable union of open interval. Consider a function $\phi$ with compact support $K_\phi\subset \mathcal I^+$. Multiplying~\eqref{seq} by $\phi$ and integrating in $K_\phi$ we obtain
\begin{equation}\label{eq33}
\int_{K_\phi} v_m'(t) \phi'(t) \,dt= \int_{K_\phi} \frac{g_{n_m}(L_m,t,x_m(t))}{x_m(t)} \, v_m(t) \phi(t) \, dt\,.
\end{equation}
By compactness we have $\min_{K_\phi} v > \hat\delta$, for a suitable $\hat\delta>0$, thus giving us that
$\lim_m x_m(t) = +\infty$ uniformly for $t\in K_\phi$. Assuming, up to subsequence, $x_m(t)>1$ for every $t\in K_\phi$ we 
can find, for every $j>0$ an index $m_j$ such that
$$
\mu_N - \frac1j <
 \frac{g_{n_{m_j}}(L_{m_j},t,x_{m_j}(t))}{x_{m_j}(t)} 
 < \mu_{N+1} + \frac 1j\,.
$$
Hence, the subsequence
$$
\left( \frac{g_{n_{m_j}}(L_{m_j},t,x_{m_j}(t))}{x_{m_j}(t)} \right)_j
$$
is bounded in $L^2(K_\phi)$, so up to a subsequence it converges weakly to a certain function $p(t)$ such that $\mu_N\leq p(t) \leq \mu_{N+1}$ for almost every $t\in K_\phi$ almost everywhere.  
Hence, passing to the limit in~\eqref{eq33}, we obtain
$$
\int_{K_\phi} v'(t) \phi'(t) \, dt = \int_{K_\phi} p(t) v(t) \phi(t) \, dt\,.
$$
It is possible to extend the function $p$ to the whole set $\mathcal I^+$, so that
$$
\int_{\mathcal I^+} v'(t) \phi'(t) \, dt = \int_{\mathcal I^+} p(t) v(t) \phi(t) \, dt\,,
$$
thus giving us that $v$ is a {\em weak} solution of $v'' + p(t) v =0$ in $\mathcal I^+$. In particular $v\in H^2_{loc}(\mathcal I^+)$ and $v\in C^1(\mathcal I^+)$. We must show that $p(t) = \mu_{N+1}$ for almost every $t\in \mathcal I^+$.

We recall that the functions $x_m$ perform in the phase plane exactly $N+1$ rotations around the origin, so there exist
$$
\alpha_1^m < \beta_1^m < \alpha_2^m < \beta_2^m < \cdots < \alpha_{N+1}^m < \beta_{N+1}^m <  \alpha_{N+2}^m = \alpha_1^m + T
$$
such that, for every $r\in\{1,\ldots,N+1\}$,
$$
x_m(t) >0 \text{ for every } t\in(\alpha_r^m, \beta_r^m)\,, $$ $$
x_m(t) <0 \text{ for every } t\in(\beta_r^m, \alpha_{r+1}^m)\,.
$$
Up to subsequences, we can assume that $\alpha_r^m \to \check\xi_r$ and $\beta_r^m \to \hat\xi_r$ such that
$$
\check\xi_1 \leq \hat\xi_1 \leq \check\xi_2 \leq \hat\xi_2 \leq \cdots \leq\check\xi_{N+1} \leq \hat\xi_{N+1} \leq \check\xi_{N+2} = \check\xi_1 +T \,.
$$
Being $\alpha_{r+1}^m - \beta_r^m < \pi/ \sqrt {n_m}$, then $\check\xi_{r+1} = \hat\xi _r $. Moreover, using the estimate in~\eqref{goodtime}, from $\beta^m_r - \alpha^m_r > \pi/ \sqrt{\mu_{N+1}} - \varepsilon = T/(N+1) - \varepsilon$, we have 
$\hat\xi_r - \hat\xi_{r-1} = T/(N+1)$.
From $v_n(\alpha^m_r)=v_m(\beta^m_r)=0$ we have $v(\xi_r)=0$, where we denote $\xi_r = \check\xi_r = \hat\xi_{r+1}$.

\medbreak

We now consider an interval $[\alpha,\beta]$ with $v(\alpha)=v(\beta)=0$ and $v(t)>0$ in $(\alpha,\beta)$. Writing in polar coordinates
$$
\begin{cases}
v(t) = \hat\rho(t) \cos( \hat\vartheta(t))\\
v'(t) = \hat\rho(t) \sin( \hat\vartheta(t)) \,,
\end{cases}
$$
we obtain the expression of the angular velocity of the solution in the phase plane
$$
- \hat\vartheta'(t) = \frac{p(t) v(t)^2 + v'(t)^2}{v(t)^2 + v'(t)^2} \,,
$$
thus obtaining
$$
\frac{- \hat\vartheta'(t)}{\mu_{N+1} \cos^2(\hat\vartheta(t)) + \sin^2(\hat\vartheta(t))} 
\leq 1\leq
\frac{- \hat\vartheta'(t)}{\mu_N \cos^2(\hat\vartheta(t)) + \sin^2(\hat\vartheta(t))}\,,
$$
so that
$$
\frac{T}{N+1} \leq \beta - \alpha \leq \frac TN.
$$
The previous reasoning can be done for every interval $[\xi_r,\xi_{r+1}]$, so the only reasonable conclusion is that $\beta-\alpha = T/(N+1)$. Hence, if the solution $v$ becomes positive, then it is forced to remain positive for a time equal to $T/(N+1)$. Passing to modified polar coordinates
$$
\begin{cases}
v(t) = \frac{1}{\sqrt{\mu_{N+1}}}\tilde\rho(t) \cos( \tilde\vartheta(t))\\
v'(t) = \tilde\rho(t) \sin( \tilde\vartheta(t)) \,,
\end{cases}
$$
we obtain integrating $-\tilde\vartheta'$ on $[\alpha,\beta]$
$$
\pi = \sqrt{\mu_{N+1}} \int_\alpha^\beta \frac{p(t)v(t)^2 + v'(t)^2}{\mu_{N+1}v(t)^2 + v'(t)^2} \,dt \leq \sqrt{\mu_{N+1}} \frac{T}{N+1} = \pi\,,
$$
thus giving us $p(t)=\mu_{N+1}$ for almost every $t\in[\alpha,\beta]$. In particular, $\mathcal I^+$ is the union of the intervals $(\xi_r,\xi_{r+1})$, where we have
\begin{equation}\label{thisform}
v(t) = c_r \sin (\sqrt{\mu_{N+1}} (t- \xi_r))\,,
\end{equation}
with $c_r\in[0,1]$ and at least one of them is equal to $1$, being $\|v\|_\infty=1$.

\medbreak

We now prove that 
\begin{equation}\label{cr1}
c_r=1 \text{ for every } r\in\{1,\ldots,N+1\} \,.
\end{equation}
\medbreak

The functions $v_m$ solve equation~\eqref{seq}, which we rewrite in the simpler form
$$
v_m'' + h_m(t,v_m) =0\,,
$$
where, for every $m$,
\begin{equation}\label{hgood}
|h_m(t,v)|\leq d (v + 1) \quad\text{ for every } t\in[0,T] \text{ and } v\geq 0 \,,
\end{equation}
for a suitable constant $d>0$.

We show that, if $v$ is positive in $[a,b]\subset[0,T]$ then $v_m$ $C^1$-converges to $v$. We have already seen that $(v_m)_m$ is bounded in $C^1$, and by~\eqref{hgood} as an immediate consequence we get $|v_m''(t)|\leq |h_m(t,v_m)|\leq 2d$ for every $t\in[a,b]$. So, being $v_m$ bounded in $C^2$ in such a interval, by the Ascoli-Arzel\` a theorem, we have that $v_m$ $C^1$-converges to $v$ in $[a,b]$.

The $C^1$-convergence and the estimate in~\eqref{hgood} are the ingredients we need to prove that the solution $v$ has only isolated zeros $\xi_r$.

We start proving that if the left derivative $-v'(\xi_r^-)=\eta>0$ for a certain index $r$, then $-v_m'(\beta_r^m)>\eta/2$ for $m$ large enough.
For every $\epsilon_0>0$, we can find $0<s_1<s_2$ sufficiently small to have
$$
\frac12 \epsilon_0  <v(\xi_r - s) < \frac32 \epsilon_0 \text{ and } |v'(\xi_r - s) + \eta| < \epsilon_0\,,
$$
for every $s\in(s_1,s_2)$. Being $v_m$ $C^1$-convergent to $v$ in $(\xi_r - s_2,\xi_r - s_1)$, for $m$ large enough,
$$
\frac12  \epsilon_0  <v_m(\xi_r - s) < \frac32 \epsilon_0 \text{ and } |v_m'(\xi_r - s) + \eta| < 2 \epsilon_0\,,
$$
for every $s\in(s_1,s_2)$. Being $|v_m''| \leq 2d$ in this interval, we find
$$
v_m'^2(\beta_r^m) \geq (\eta -2\epsilon_0)^2 - 6d\epsilon_0 > \eta^2/4
$$
choosing $\epsilon_0$ sufficiently small.

We prove now that if $-v'(\xi_r^-)=\eta>0$ for a certain index $r$, then $\xi_r$ is an isolated zero of $v$.
Suppose by contradiction that there exists $\ee_0\in(0,\eta/8d)$, with $d$ as in~\eqref{hgood}, such that $v(\xi_r+\ee_0)=0$. For every $m$ large enough we have 
$|\alpha_r^m-\xi_r|<\ee_0/4$ and by the previous computation $v_m'(\alpha_r^m)=-v_m'(\beta_r^m)>\eta/2$.

The property that $v_m'(\alpha_r^m)=-v_m'(\beta_r^m)$ follows directly by the fact that the nonlinearities $g_n$ do not depend by $t$ when $x<0$, leads us to treat a Landesman-Lazer condition involving function $\psi_j$ as in~\eqref{psigood} rather than in~\eqref{psibad}. Here lies one of the main differences between our result and the one obtained by Fonda and Garrione in~\cite{FG2}.

Being $|v_m''|\leq 2d$ when $v_m$ is positive, we can show that if $s<\eta/4d$ then $v_m(\alpha_r^m+s)>s\, \eta/4$. By construction $\xi_r+\ee_0=\alpha_r^m+s_0$ for a certain $s_0\in(\ee_0/2,\eta/4d)$, so that we obtain $v_m(\xi_r+\ee_0)=v_m(\alpha_r^m+s_0)>\eta\ee_0/8$ for every $m$ large enough, thus contradicting $v_m\to v$.

\medbreak

 Being $v$ as in~\eqref{thisform}, and $c_r=1$ for at least one value $r\in\{1,\ldots,N+1\}$, for such index 
 $v'(\xi_{r+1}^-) < 0$ holds. The previous reasoning gives us that $c_{r+1}>0$. Iterating the procedure we can prove that $c_r>0$ for every index $r\in\{1,\ldots,N+1\}$.

We prove now that for every $r\in\{1,\ldots,N+1\}$, the left and right derivatives satisfy $v'(\xi_{r}^-)=-v'(\xi_{r}^+)$, thus we can conclude that $c_r=1$ for every $r\in\{1,\ldots,N+1\}$.

\medbreak

Suppose by contradiction that there exists $r\in\{1,\ldots,N+1\}$ such that
$$
v'(\xi_{r}^+) + v'(\xi_{r}^-) \neq 0\,.
$$
Without loss of generality we suppose this value to be positive. The other case follows similarly. So, assume
$$
v'(\xi_{r}^+) + v'(\xi_{r}^-) > p_0  > 0 \text{ and } 0<q_0 \leq \min\{ -v'(\xi_{r}^-), v'(\xi_{r}^+)\}\,.
$$
Arguing as above, for every $\epsilon_0>0$, we can find $0<s_1<s_2$ sufficiently small to have
$$
\frac12 \epsilon_0  <v(\xi_r \pm s) < \frac32 \epsilon_0 \text{ and } |v'(\xi_r \pm s) - v'(\xi_r^\pm)| < \epsilon_0\,,
$$
for every $s\in(s_1,s_2)$. Being $v_m$ $C^1$-convergent to $v$ in $(\xi_r - s_2,\xi_r - s_1)$ and in $(\xi_r + s_1,\xi_r + s_2)$, for $m$ large enough,
$$
\frac12  \epsilon_0  <v_m(\xi_r \pm s) < \frac32 \epsilon_0 \text{ and } |v_m'(\xi_r \pm s) - v'(\xi_r^\pm)| < 2 \epsilon_0\,,
$$
for every $s\in(s_1,s_2)$. Being $|v_m''| \leq 2d$, we find
$$
[ v'(\xi_r^-) - 2\epsilon_0]^2 + 4\epsilon_0 d \geq v_m'^2(\beta_r^m)= v_m'^2(\alpha_{r+1}^m) \geq [  v'(\xi_r^+) - 2\epsilon_0]^2 - 4\epsilon_0 d\,.
$$
Hence
$$
0<p_0 < v'(\xi_{r}^+) + v'(\xi_{r}^-) \leq \frac{4 \epsilon_0 d}{q_0} + 4\epsilon_0\,,
$$
thus giving us a contradiction for $\epsilon_0$ sufficiently small.

\medbreak

We have proved that $v$, in every interval $[\alpha,\beta]=[\xi_r,\xi_{r+1}]$, satisfies
\begin{equation}\label{howisv}
v(t) =  \sin \left( \sqrt{\mu_{N+1}}\,(t-\alpha)\right).
\end{equation}
So, $v$ is a solution of the following Dirichlet problem:
$$
\begin{cases}
v'' + \mu_{N+1} \, v =0\\
v(\alpha)=0, \quad v(\beta)=0\,.
\end{cases}
$$
Let us consider the orthonormal basis $(\phi_k)_k$ of $L^2(\alpha,\beta)$ made of the eigenfunctions solving the Dirichlet problem
$$
\begin{cases}
\phi_k'' + \mu_{k} \, \phi_k =0\\
\phi_k(\alpha)=0, \quad \phi_k(\beta)=0\,,
\end{cases}
$$
where $\mu_k=(k \pi/T)^2$ is the $k$-th eigenvalue. Denoting by $\pscal \cdot\cdot$ and $\| \cdot \|_2$, respectively, the scalar product and the norm in $L^2(\alpha,\beta)$, we can write the Fourier series of all the functions $x_n$ and split it as follows:
\beq
x_m &=& \sum_{k=1}^{\infty} \pscal{x_m}{\phi_k} \phi_k \\
&=&  \underbrace{\pscal{x_m}{\phi_{N+1}}\phi_{N+1} }_{x_m^0}  + \underbrace{ \sum_{k\neq N+1} \pscal{x_m}{\phi_k} \phi_k }_{x_m^\perp}\,,
\eeq
with the following property
$$
(x_m'')^0 = (x_m^0)'' \quad \text{ and } \quad (x_m'')^\perp = (x_m^\perp)'' .
$$
Moreover, one has $v_m= v_m^0 + v_m^\perp$ with
$$
v_m^0 = \frac{x_m^0}{\|x_m \|_\infty} \quad\text{ and } \quad v_m^\perp = \frac{x_m^\perp}{\|x_m \|_\infty}\,.
$$
Being $v=\|v\|_2 \phi_{N+1}$, we have $v_m^0 \to v$ uniformly in $[\alpha,\beta]$. Moreover $v_m^0\geq 0$, for $m$ sufficiently large.

Multiplying equation~\eqref{gnm} by $v_m^0$ and integrating in the interval $[\alpha,\beta]$ we obtain
\begin{equation}\label{1est}
\begin{array}{l}
\ds \int_\alpha^\beta g_{n_m}(L_m,t,x_m(t)) v_m^0(t) \,dt 
	= -\int_\alpha^\beta (x_m^0)''(t) v^0_m(t)\,dt\\ [3mm]
\ds\hspace{10mm}=  - \int_\alpha^\beta x_m^0(t) (v_m^0)''(t) \,dt
=   \int_\alpha^\beta \mu_{N+1} x_m^0(t) v_m^0(t) \,dt \\ [3mm]
 \ds\hspace{20mm}=  \int_\alpha^\beta \mu_{N+1} x_m(t) v_m^0(t) \,dt \,.
 \end{array}
\end{equation}
Defining $r_m(t,x)= g_{n_m}(L_m,t,x) - \mu_{N+1} x$ we have
$$
\int_{\alpha}^{\beta} r_m(t,x_m(t)) v_m^0(t) \,dt = 0\,,
$$
and applying Fatou's lemma
$$
\int_{\alpha}^{\beta} \limsup_{m\to\infty} r_m(t,x_m(t)) v_m^0(t) \,dt \geq 0\,.
$$
It is easy to see that for every $s_0\in(\alpha,\beta)$ it is possible to find $m(s_0)$ such that $x_m(s_0)>1$ for every $m>m(s_0)$.
Hence, being $v_m^0 \to v$ and $L_m \to L_\dag$,
we have
\begin{equation}\label{LLnot}
\int_{\alpha}^{\beta} \limsup_{x\to+\infty} [f(t,x+R_0)-\mu_{N+1} x] v(t) \,dt \geq 0\,.
\end{equation}
The previous estimate can be obtained for every interval $(\xi_r,\xi_{r+1})$, thus contradicting~\eqref{LLcond} setting $\tau=\xi_1$.

\medbreak

We have so finished to prove the case in which the sequence $(x_m)_m$ consists of solutions performing $N+1$ rotations in the phase-plane around the origin for an infinite number of index $m$. We now treat the case in which the solutions perform $N$ rotations around the origin.

\medbreak

In this case we have to prove that the limit function $v$ solves $v'' + \mu_N v =0$ for almost every $t$. With the same procedure we can show that there exists a $L^2$-function $q(t)$, satisfying $\mu_N \leq q(t)\leq \mu_{N+1}$ almost everywhere, such that $v$ is a weak solution of $v'' + q(t)v =0$. Then, with a similar procedure, it is possible to introduce some instants $\alpha_r^m$, $\beta_r^m$, when $v_m$ vanishes, converging to some values $\xi_r$ (with $r\in\{1,\ldots,N\}$). Unfortunately. we cannot conclude immediately that  $\xi_{r+1}-\xi_r = T/N$ holds. In fact, using~\eqref{goodtime}, we can verify only that $\xi_{r+1}-\xi_r \leq T/N$. Nevertheless, we have that $v_m$ $C^1$-converges to $v$ when $v$ is positive. Moreover, following the reasoning which gave us the estimates in~\eqref{goodderivative}, we can prove that for at least one index $r\in\{1,\ldots,N\}$ the left-derivative satisfies $v'(\xi_r^-)<0$ (e.g. the index in which $v$ attains the maximum). Arguing as above, we can prove that whenever $v'(\xi_r^-)<0$ for a certain index $r\in\{1,\ldots,N\}$ then $\xi_r$ is an isolated zero and $v'(\xi_r^+)=-v'(\xi_r^-)>0$. Iterating this reasoning, we have that all the $\xi_r$ are isolated zeros of $v$ thus obtaining the needed estimate $\xi_{r+1}-\xi_r =T/N$. Now, with a similar procedure, we can prove that $p(t)=\mu_N$, thus giving us that $v(t)= c_r \sin(\sqrt{\mu_N}(t-\xi_r))$ in the interval $(\xi_r,\xi_{r+1})$. Also in this case, in the same way, it is possible to conclude that $c_r=1$ for every $r\in\{1,\ldots,N\}$. Then, we can consider an interval $[\alpha,\beta]$ with $v>0$ in $(\alpha,\beta)$, and with a similar reasoning we can obtain a {\em liminf} estimate similar to the one obtained~\eqref{LLnot}, thus gaining a contradiction with assumption~\eqref{LLcond2}.

The proposition is thus proved. \cvd


\medbreak

Let us spend few words about the possibility of extending Theorem~\ref{LLtheorem} to the case $N=0$, where $\mu_0=0$ and $\mu_1=(\pi/T)^2$. 
There is a wide literature (cf.~\cite{FZ,F,GO,HOZ,MawWar}) treating nonlinearities {\em lying} under the first curve of the Dancer-Fu\v c\'{\i}k spectrum. The Landesman-Lazer condition~\eqref{LLcond2} in this case reduces to a sign condition on the nonlinearity $f$. Unfortunately, it is not possible to obtain a proof with the same procedure, cf. the estimate in~\eqref{rt}. For briefness we do not enter in such details in this paper. However, let us state the following weaker result for a nonlinearity with a {\em one-side} resonance condition, the proof of which works similarly to the one of Theorem~\ref{LLtheorem}.

\begin{theorem}\label{LLtheoremN0}
Assume that there exists a constant $\bar\ee$ such that
$$
\liminf_{x\to \infty} \frac{f(t,x)}{x} \geq \bar\ee > 0\,,
$$
uniformly for every $t\in[0,T]$ and that
there exists a constant $\hat\eta$ such that, for $N>0$,
$$
f(t,x)\leq \mu_1 x +\hat\eta
$$
for every $t\in[0,T]$ and every $x>R_0$. Moreover, for every $\tau\in[0,T]$,
\begin{equation}\label{LLcond2N0}
\int_0^T \liminf_{x\to +\infty} ( f(t,x+R_0) - \mu_1 x ) \psi_{1}(t+\tau) \, dt > 0\,.
\end{equation}
Then, for every integer $\nu>0$, there exists an integer $k_\nu>0$, such that for every integer $k\geq k_\nu$ there exists at least one periodic $R_0$-bouncing solution $(\rho,\vartheta)$ of~\eqref{sist1} with period $kT$, which makes exactly $\nu$ revolutions around the origin in the period time $kT$, i.e. satisfying~\eqref{percond}.
\end{theorem}

\section{Systems on cylinders}\label{seccil}

In this section we briefly explain how the previous results could be applied to a class of systems defined in $\R^{d_1+d_2}$, modeling a  particle hitting a cylinder $\mathbb S^{d_1-1}\times \R^{d_2}$. 
The case $d_1=2$ and $d_2=1$ models bounces on a proper cylinder.
For briefness we will present the result for {\em non-resonant} nonlinearities. We consider the differential equations
\begin{equation}\label{eqcil}
\begin{array}{l}
\ds {\bf x}'' + {\text{f}_1}(t,|{\bf x}|) {\bf x}+{\text{\bf b}_1}(t,{\bf x},{\bf y})=0\,,\\
\ds {\bf y}'' + {\text{\bf f}_2}(t,{\bf y}) +{\text{\bf b}_2}(t,{\bf x},{\bf y})=0\,,
\end{array}
\end{equation}
where ${\bf x}\in \mathbb S^{d_1-1}$ and ${\bf y}\in\R^{d_2}$.
We assume for simplicity all the functions to be continuous.
We suppose that $f=\text{f}_1(t,\rho) \rho$ satisfies the assumptions of Theorem~\ref{main}.
The function ${\text{\bf b}_1}:\R^{1+d_1+d_2} \to \R^{d_1}$ satisfies ${\text{\bf b}_1}(t,{\bf x},{\bf y})=b_1(t,{\bf x},{\bf y}) {\bf x}$ with
$$
\lim_{|{\bf x}|\to \infty} \frac{{b_1}(t,{\bf x},{\bf y}) }{|{\bf x}|} = 0
$$
uniformly in $t$ and {\bf y}. Assume that the second equation in~\eqref{eqcil} can be viewed, in every components, as
$$
y_i'' + f_{2,i}(t,y_i) +b_{2,i}(t,{\bf x},{\bf y}) = 0\,,
$$
where
$$
\lim_{|y_i |\to \infty} \frac{{b_{2,i}}(t,{\bf x},{\bf y}) }{y_i} = 0
$$
uniformly in all the other variables. We can assume, as an example of application,

$$
\check \mu_i \leq \liminf_{y_i\to+\infty} \frac{f_{2,i}(t,y_i)}{y_i} \leq \limsup_{y_i \to+\infty} \frac{f_{2,i}(t,y_i)}{y_i} \leq \hat\mu_i
$$
$$
\check \nu_i \leq \liminf_{y_i\to-\infty} \frac{f_{2,i}(t,y_i)}{y_i} \leq \limsup_{y_i \to-\infty} \frac{f_{2,i}(t,y_i)}{y_i} \leq \hat\nu_i
$$
uniformly in $t$,
with
$$
\frac{T}{(N_i+1)\pi} <\frac{1}{\sqrt{\hat\mu_i}}+\frac{1}{\sqrt{\hat\nu_i}} \leq \frac{1}{\sqrt{\check\mu_i}}+\frac{1}{\sqrt{\check\nu_i}}< \frac{T}{N_i\pi}
$$
for some positive constant $\check\mu_i$, $\check\nu_i$, $\hat\mu_i$, $\hat\nu_i$ and an integer $N_i>0$.

\begin{theorem}
Under the previous assumptions, for every integer $\ell>0$, there exists an integer $k_\ell>0$, such that for every integer $k\geq k_\ell$ there exists at least one periodic solution $(\text{\bf x},\text{\bf y})$ of~\eqref{eqcil}, such that {\bf x}, can be parametrized in polar coordinates $(\rho,\vartheta)$ and $\rho$ is a $R_0$-bouncing solution. Such solutions satisfies the following periodicity conditions:
$$
\begin{array}{l}
\ds \rho(t+T)=\rho(t)\,,\\
\ds \vartheta(t+kT)=\vartheta(t)+2\pi\nu\,,\\ 
\ds \text{\bf y}(t+T)=\text{\bf y}(t)\,.
\end{array}
$$
\end{theorem}

The proof of such a result can be obtained by glueing together the results contained in this paper (for the $\text{\bf x}$ coordinate) and classical results (for the $\text{\bf y}$ coordinate). The key tool is the fact that the equations are {\em weakly} coupled.

\newcommand\mybib[8]{{\bibitem{#1} {\sc #2}, {\em #3}, {#4}~{\bf #5} ({#6}), {#7}--{#8}.}} 
\newcommand\mybibb[4]{{\bibitem{#1} {\sc #2}, {\em #3}, {#4}}} 

\noindent Author's address:

\bigbreak

\begin{tabular}{l}
Andrea Sfecci\\
Universit\`a Politecnica delle Marche\\
Dipartimento di Ingegneria Industriale e Scienze Matematiche\\
Via Brecce Bianche 12\\
60131 Ancona\\
Italy\\
e-mail: sfecci@dipmat.univpm.it
\end{tabular}

\medbreak

\noindent Mathematics Subject Classification: 34B15, 34C25

\medbreak
\noindent Keywords: Keplerian problem, periodic solutions, impact, resonance,\\ Landesman-Lazer condition.

\end{document}